\documentclass[11pt]{amsart}
\usepackage{amsmath}
\usepackage{graphicx}
\numberwithin{equation}{section}
\newtheorem{theorem}{Theorem}
\newtheorem{lemma}[theorem]{Lemma}
\newtheorem{remark}{Remark}
\newtheorem{corollary}[theorem]{Corollary}

{\rm}
{\rm}
{\rm}

\def\N{\mathbb{N}}
\def\R{\mathbb{R}}
\def\K{\mathbf{K}}

\def\va{\vert\alpha\vert}
\begin{document}

\title{Sufficient conditions  for a real polynomial to be a sum of squares}

\author{Jean B. Lasserre}
\address{LAAS-CNRS and Institute of Mathematics,
 LAAS, 7 avenue du Colonel Roche,
 31077 Toulouse cedex 4, France}
\email{lasserre@laas.fr}
\date{}
\begin{abstract}
We provide explicit sufficient conditions for a polynomial $f$ 
to be a sum of squares (s.o.s.), linear in the coefficients of $f$.
All conditions are simple and provide an explicit description of
a convex polyhedral subcone of the cone of s.o.s. polynomials 
of degree at most $2d$. We also provide a simple condition to ensure that $f$ is s.o.s., possibly after adding a constant.
\end{abstract}

\maketitle
\section{Introduction}
The cone $\Sigma^2\subset\R[X]$ of real polynomials that are 
sum of squares (s.o.s.) and its subcone $\Sigma^2_d$
of s.o.s. of degree at most $2d$,
play a fundamental role
in many areas, and particularly in optimization; see for instance
Lasserre \cite{lasserre1,lasserre2}, Parrilo \cite{parrilo} and Schweighofer
\cite{markus}.
When considered as a convex cone of a finite dimensional euclidean space,
$\Sigma^2_d$ has a {\it lifted semidefinite representation} (such sets are called 
SDr sets in \cite{bental}). 
That is, $\Sigma^2_d$ is the projection of a convex cone of
an euclidean space of higher dimension, defined in terms of the coefficients
of the polynomial and additional variables (the "lifting"). However,
so far there is no simple description of $\Sigma^2_d$ given {\it directly} in terms of 
the coefficients of the polynomial. 
For more details on SDr sets, the interested reader is referred to e.g. Ben Tal and Nemirovski \cite{bental}, Helton and Vinnikov \cite{helton},
Lewis et al. \cite{lewis}.

Of course, one could use Tarski's quantifier elimination to
provide a description of $\Sigma^2_d$,
solely in terms of the coefficients,
but such a description is likely hopeless to be {\it simple}; in particular,
it could be sensitive to the degree $d$.
Therefore,  a more reasonable goal is to
search for simple descriptions of {\it subsets} (or {\it subcones}) of $\Sigma_d^2$ only.
This is the purpose of this note in which we provide simple sufficient conditions for a polynomial $f\in\R[X]$ of degree at most $2d$, to be s.o.s. All conditions are expressed directly in terms 
of the coefficients $(f_\alpha)$, with no additional variable
(i.e. with no lifting) and define a convex polyhedral subcone of $\Sigma^2_d$.
Finally, we also provide a sufficient condition on the coefficients of highest degree
to ensure that $f$ is s.o.s., possibly after adding a constant.
All conditions stress the importance
of the {\it essential}
monomials $(X_i^{2k})$ which also play an important role
for approximating nonnegative polynomials by s.o.s.,
as demonstrated in e.g. \cite{lasserre2,lasserre-tim}.

\section{Conditions for being s.o.s.}
For $\alpha\in\N^n$ let $\va:=\sum_{i=1}^n\vert\alpha_i\vert$.
Let $\R[X]$ be the ring of real polynomials in the variables
$X=(X_1,\ldots,X_n)$, and let $\R[X]_{2d}$ the vector space of 
real polynomials of degree at most $2d$, with canonical basis of monomials
$(X^\alpha)=\{X^\alpha\::\:\alpha\in\N^n;\: \va\leq 2d\}$.
Given a sequence $y=(y_\alpha)\subset\R$ indexed
in the canonical basis $(X^\alpha)$, let $L_y:\R[X]_{2d}\to\,\R$ be the linear
mapping
\[f\,(=\sum_\alpha f_\alpha\,X^\alpha)\quad\mapsto\quad L_y(f)\,=\,\sum_{\alpha}f_\alpha\,y_\alpha,\quad f\in\R[X]_{2d},\]
and let $M_d(y)$ be the {\it moment} matrix with rows and columns indexed in $(X^\alpha)$, and defined by
\begin{equation}
\label{moment}
M_d(y)(\alpha,\beta)\,:=\,L_y(X^{\alpha+\beta})\,=\,y_{\alpha+\beta},\qquad\alpha,\beta\in\N^n:\: \vert\alpha\vert,\vert\beta\vert\,\leq\,d.
\end{equation}
Let the notation $M_d(y)\succeq0$ stand for $M_d(y)$ is positive semidefinite.
It is clear that 
\[M_d(y)\,\succeq\, 0\quad\Longleftrightarrow\quad L_y(f^2)\geq\,0\quad\forall\,f\in\R[X]_d.\]
The set $\Sigma^2_d\subset\R[X]_{2d}$ of s.o.s. polynomials of degree at most $2d$
is a finite-dimensional convex cone, and 
\begin{equation}
\label{sos}
f\in\Sigma^2_d\quad\Longleftrightarrow\quad L_y(f)\,\geq\,0\quad
\forall \,y\:\mbox{ s.t. }M_d(y)\,\succeq\,0.
\end{equation}
\begin{remark}
\label{rem1}
{\rm To prove that $L_y(f)\geq0$ for all $y$ such that $M_d(y)\succeq0$ 
it suffices to prove that $L_y(f)\geq0$ for all $y$ such that 
$M_d(y)\succeq0$ and $L_y(1)>0$ (and equivalently, by homogeneity,
for all $y$ such that $M_d(y)\succeq0$ and $L_y(1)=1$).

Indeed, suppose that $L_y(f)\geq0$ for all $y$ such that
$M_d(y)\succeq0$ and $L_y(1)>0$. Next, let $y$ be such that $M_d(y)\succeq0$ and $L_y(1)=0$. Fix $\epsilon>0$ arbitrary and let $y(\epsilon):=y+(\epsilon,0,\ldots,0)$ so that $L_{y(\epsilon)}(X^\alpha)=y_\alpha$ if $\alpha\neq0$ and $L_{y(\epsilon)}(1)=\epsilon >0$. Therefore $M_d(y(\epsilon))\succeq0$ (because $M_d(y)\succeq0$) and so 
$0\leq L_{y(\epsilon)}(f)=\epsilon f_0+L_y(f)$. As $\epsilon>0$ was arbitrary, letting $\epsilon\downarrow 0$ yields the desired result
$L_y(f)\geq0$.
}\end{remark}
We first recall a preliminary result whose proof can be found in 
Lasserre and Netzer \cite{lasserre-tim}.
\begin{lemma}[\cite{lasserre-tim}]
\label{lemma1}
With $d\geq1$, let $y=(y_\alpha)\subset\R$ be such that the moment matrix 
$M_d(y)$ defined in (\ref{moment}) is positive semidefinite, and let
$\tau_d:=\displaystyle\max_{i=1,\ldots,n}\,L_y(X_i^{2d})$.  Then:
\begin{equation}
\label{lemma1-1}
\vert L_y(X^\alpha)\vert\,\leq\,\max[\,L_y(1)\,,\,\tau_d\,],\qquad\forall\,\alpha\in\N^n:
\quad \va\leq 2d.
\end{equation}
\end{lemma}

We next complement Lemma \ref{lemma1}.

\begin{lemma}
\label{lemma2}
Let $y=(y_\alpha)\subset\R$ be normalized with $y_0=L_y(1)=1$, and such that 
$M_d(y)\succeq0$. Let
$\tau_d:=\displaystyle\max_{i=1,\ldots,n}\,L_y(X_i^{2d})$.  Then:
\begin{equation}
\label{lemma2-1}
\vert L_y(X^\alpha)\vert^{1/\vert\alpha\vert}\,\leq\,\tau_d^{1/2d},\qquad \forall\,\alpha\in\N^n:\:1\leq\vert\alpha\vert\,\leq\,2d.
\end{equation}
\end{lemma}
For a proof see \S\ref{proof-lemma2}.
\subsection{Conditions for a polynomial to be s.o.s.}
With $d\in\N$, let $\Gamma\subset\N^n$ be the set defined by:
\begin{equation}
\label{setgamma}
\Gamma\,:=\,\{\:\alpha\in\N^n\::\:\va\leq 2d;\quad
\alpha=2\beta\quad\mbox{for some }\beta\in\N^n\}.
\end{equation}
We now provide our first condition.
\begin{theorem}
\label{th1}
Let $f\in\R[X]_{2d}$ and write $f$ in the form
\begin{equation}
\label{form}
 f\,=\,f_0+\sum_{i=1}^nf_{i2d}\,X_i^{2d}\:+h,
\end{equation}
where $h\in\R[X]_{2d}$ contains no essential monomial $X_i^{2d}$. If

\begin{eqnarray}
\label{th1-1}
f_0&\geq&\sum_{\alpha\not\in\Gamma}\,\vert f_\alpha\vert \quad-
\sum_{\alpha\in\Gamma}\min[0,f_\alpha]\\
\label{th1-2}
\min_{i=1,\ldots,n}\,f_{i2d}&\geq&
\sum_{\alpha\not\in\Gamma}\vert f_\alpha\vert\,\frac{\vert\alpha\vert}{2d}\quad-
\sum_{\alpha\in\Gamma}\min[0,f_\alpha]\,\frac{\vert\alpha\vert}{2d}
\end{eqnarray}
then $f\in\Sigma^2_d$.
\end{theorem} 
For a proof see \S\ref{proof-th1}.
The sufficient conditions (\ref{th1-1})-(\ref{th1-2}) define a
polyhedral convex cone in the euclidean space of coefficients $(f_\alpha)$ of polynomials $f\in\R[X]_{2d}$. This is because the functions,
\[f\mapsto\displaystyle\min_{i=1,\ldots,n}f_{i2d}\,,\quad
f\mapsto\min[0,f_\alpha]\,,\quad
f\mapsto-\vert f_\alpha\vert,\]
are all piecewise linear and concave. The description
(\ref{th1-1})-(\ref{th1-2}) of this convex polyhedral cone is {\it explicit}
and given only in terms of the coefficients $(f_\alpha)$, i.e., with no lifting.\\

Notice that (\ref{th1-1})-(\ref{th1-2}) together with $f_{i2d}=0$ for some $i$, implies $f_\alpha=0$ for all $\alpha\not\in\Gamma$, and
$f_\alpha\geq0$ for all $\alpha\in\Gamma$, in which case $f$ is obviously s.o.s.\\

Theorem \ref{th1} is interesting when $f$ has a few non zero coefficients. When $f$ has a lot of non zero coefficients and contains the essential monomials
$X_i^{2k}$ for all $k=1,\ldots,d$, all with positive coefficients, one
provides the following alternative sufficient condition.
With $k\leq d$, let 
\begin{eqnarray}
\label{setgammak1}
\Gamma^1_k&:=&\{\:\alpha\in\N^n\::\quad 2k-1\,\leq\,\va\,\leq\,2k\:\}\\
\label{setgammak2}
\Gamma^2_k&:=&\{\:\alpha\in\,\Gamma_k^1\::\quad
\alpha=2\beta\quad\mbox{for some }\beta\in\N^n\}.
\end{eqnarray}

\begin{corollary}
\label{coro1}
Let $f\in\R[X]_{2d}$ and write $f$ in the form
\begin{equation}
\label{newform}
f\,=\,f_0+h+\sum_{k=1}^d\sum_{i=1}^n\,f_{i2k}\,X_i^{2k},
\end{equation}
where $h\in\R[X]_{2d}$ contains no essential monomial $X_i^{2k}$. If 
\begin{eqnarray}
\label{coro1-1}
\frac{f_0}{d}&\geq&\sum_{\alpha\in\Gamma^1_k\setminus\Gamma^2_k}\,\vert f_{\alpha}\vert \quad -
\sum_{\alpha\in\Gamma^2_k}\min[0,f_{\alpha}]\\
\label{coro1-2}
\min_{i=1,\ldots,n}\,f_{i2k}&\geq&
\sum_{\alpha\in\Gamma^1_k\setminus\Gamma^2_k}\vert f_\alpha\vert\,\frac{\vert\alpha\vert}{2k}\quad-
\sum_{\alpha\in\Gamma^2_k}\min[0,f_\alpha]\,\frac{\vert\alpha\vert}{2k}
\end{eqnarray}
for all $k=1,\ldots,d$, then $f\in\,\Sigma^2_d$.
\end{corollary}
For a proof see \S\ref{proof-coro1}.
Notice that (\ref{coro1-1})-(\ref{coro1-2}) together with $f_{i2k}=0$ for some $i$ and some $k\in\{1,\ldots d\}$, implies 
$f_\alpha=0$ for all $\alpha\in\Gamma^1_k\setminus\Gamma^2_k$,
and $f_\alpha\geq0$ for all $\alpha\in\Gamma_2^k$.

Several variants of Corollary \ref{coro1} can be 
derived; for instance, any other way to distribute
the constant term $f_0$ as $\sum_{k=1}^d f_{0k}$ with
$f_{0k}\neq f_0/d$, is valid and also 
provides another set of sufficient conditions. 
Consider also the case when $f$ can be written as
\[f\,=\,f_0+h+\sum_{k\in\K}\,\sum_{i=1}^n\,f_{i2k}\,X_i^{2k},\]
where $\K:=\{k\,\in\,\{1,\ldots,d\}\::\:\min_{i=1,\ldots,n}f_{i2k}\,>0\}$,
$d\in\K$, and $h\in\R[X]_{2d}$ contains no essential 
monomial $X_i^{2k}$, $k\in\K$. Then one may easily 
derive a set of sufficient conditions in the spirit of Corollary \ref{coro1}.\\

Finally, one provides a simple condition for 
a polynomial to be s.o.s., possibly after adding a constant.

\begin{corollary}
\label{coro2}
Let $f\in\R[X]_{2d}$ and write $f$ in the form
\begin{equation}
\label{newform2}
f\,=\,f_0+h+\sum_{i=1}^n\,f_{i2d}\,X_i^{2d},
\end{equation}
where $h\in\R[X]$ contains no essential monomial $X_i^{2d}$. If 
\begin{equation}
\label{coro2-1}
\min_{i=1,\ldots,n}\,f_{i2d}\,>\,\sum_{\alpha\not\in\Gamma;\:\vert\alpha\vert=2d}
\vert f_\alpha\vert-
\sum_{\alpha\in\Gamma;\:\vert\alpha\vert=2d}\min[0,f_\alpha]
\end{equation}
with $\Gamma$ as in (\ref{setgamma}), then $f+M\in\Sigma^2_d$ for some $M\geq0$.
\end{corollary}
\begin{proof}
Let $-M:=\min\,[0,\,\inf_y \,\{L_y(f)\,:\, M_d(y)\succeq0;\,L_y(1)=1\}]$. 
We prove that $M<+\infty$. Assume that $M=+\infty$, and let $y^j$ be a minimizing sequence.
One must have $\tau_{jd}:=\max_{i=1,\ldots,n} L_{y^j}(X_i^{2d})\to\infty$, as 
$j\to\infty$, otherwise if $\tau_{jd}$ is bounded by, say $\rho$,  by Lemma \ref{lemma1} one would have $\vert L_{y^j}(X^\alpha)\vert\leq \max[1,\rho]$ for all $\va\leq 2d$, and so $L_{y^j}(f)$ would be bounded, in contradiction 
with $L_{y^j}(f)\to-\infty$.
But then from Lemma \ref{lemma2}, for sufficiently large $j$, one obtains the contradiction
\begin{eqnarray*}
0>\frac{L_{y^j}(f)}{\tau_{jd}}&\geq &
\min_{i=1,\ldots,n}\,f_{i2d}\quad-\sum_{\alpha\not\in\Gamma;\:\vert\alpha\vert=2d}\vert f_\alpha\vert\quad +\sum_{\alpha\in\Gamma;\:\vert\alpha\vert=2d}\min[0,f_\alpha]\\
&&-\sum_{0\leq\va<2d}\vert f_\alpha \vert\,\tau_{jd}^{(\va -2d)/2d}\,\geq0,
\end{eqnarray*}
where the last inequality follows from (\ref{coro2-1}) and 
$\tau_{jd}^{(\va -2d)/2d}\to 0$ as $j\to\infty$.

Hence, $M<+\infty$ and so $L_y(f+M)\geq0$ for every $y$ with $M_d(y)\succeq0$, and $L_y(1)=1$ .
But then, in view of Remark \ref{rem1}, $L_y(f+M)\geq0$ for all
$y$ such that $M_d(y)\succeq0$, which in turn implies that
$f+M$ is s.o.s.
\end{proof}
In Theorem \ref{th1}, Corollary \ref{coro1} and \ref{coro2},
it is worth noticing the crucial role played by the 
constant term and the essential monomials  $(X_i^\alpha)$, as was already the case
in \cite{lasserre2,lasserre-tim} for approximating nonnegative polynomials by s.o.s.

 \section{Proofs}

The proof of Lemma \ref{lemma2} first requires the following auxiliary result.

\begin{lemma} \label{lemma3}
Let $d\geq1$, and $y=(y_\alpha)\subset\R$ be such that the moment matrix 
$M_d(y)$ defined in (\ref{moment}) is positive semidefinite, and let
$\tau_d:=\displaystyle\max_{i=1,\ldots,n}\,L_y(X_i^{2d})$.  Then:
$L_y(X^{2\alpha})\,\leq\,\tau_d$ for all $\alpha\in\N^n$ with $\vert\alpha\vert=d$.
 \end{lemma}
 \begin{proof} 
 The proof is by induction on the number \(n\) of variables.
The case \(n=1\) is trivial and the case 
\(n=2\) is proved in Lasserre and Netzer \cite[Lemma 4.2]{lasserre-tim}.

Let the claim be true for $k=1,\ldots,n-1$ and consider the case $n>2$.
   By the induction hypothesis, the
   claim is true for all \(L_y\left(X^{2\alpha}\right)\), where \(|\alpha| = d\)  and \(\alpha_{i}
   =0\) for some $i$. Indeed, \(L_y\) restricts to a linear form on the ring of polynomials with \(n-1\) indeterminates and satisfies all the assumptions needed. So the induction hypothesis gives the boundedness of all those values \(L_y\left(X^{2\alpha}\right)\).

 Now take \(L_y\left(X^{2\alpha}\right)\), where \(|\alpha|=d \) and all \(\alpha_{i} \geq 1\). With no loss of generality, assume \(\alpha_{1}\leq \alpha_{2} \leq...\leq\alpha_{n}.\)
   Consider the two elements \[\gamma:=(2\alpha_{1},0,\alpha_{3}+\alpha_{2}-\alpha_{1},\alpha_{4},...,\alpha_{n})\in\N^{n} \mbox { and }\]
   \[\gamma^{'}:=(0,2\alpha_{2},\alpha_{3}+\alpha_{1}-\alpha_{2},\alpha_{4},...,\alpha_{n})\in\N^{n}.\]
   We have \(|\gamma|=|\gamma^{'}|=d\) and
   \(\gamma_{2}=\gamma^{'}_{1}=0\), and from what precedes, 
   \[L_y(X^{2\gamma})\leq\tau_d\mbox{ and }
   L_y(X^{2\gamma^{'}})\leq\tau_d.\] As $M_d(y)\succeq0$,
      one also has \[L_y(X^{2\alpha})^{2}= L_y(X^{\gamma +\gamma^{'}})^{2}\leq
   L_y(X^{2\gamma})\cdot L_y(X^{2\gamma^{'}}) \leq\tau_d^{2},\] which yields the desired result $\vert L_y(X^{2\alpha})\vert\leq\tau_d$. \end{proof}
   
\subsection{Proof of Lemma \ref{lemma2}} 
\label{proof-lemma2}    
The proof is by induction on $d$. Assume it is true for $k=1,\ldots,d$, and write $M_{d+1}(y)$ in the following block form below
with appropriate matrices $V,U_i,V_i,S_i$:
\begin{center}
\begin{table}[ht!]
\begin{tabular}{|c|c|c|c|}
\hline
&&&\\
$M_{d-2}(y)$ & $U_1$ & $U_2$ & $V$\\
&&&\\
\hline
&&&\\
$U_1^T$& $S_{2d-2}$ & $V_{2d-1}$ & $V_{2d}$\\
&&&\\
\hline
&&&\\
$U_2^T$ & $V_{2d-1}^T$ & $S_{2d}$ & $V_{2d+1}$\\
&&&\\
\hline
&&&\\
$V^T$& $V_{2d}^T$ & $V_{2d+1}^T $ & $S_{2d+2}$\\
&&&\\
\hline
\end{tabular}
\end{table}
\end{center}
When $d=1$, the blocks $M_{d-2}(y)$, and $U_1,U_2,U_1^T,U_2^T,V$ disappear.\\

$\bullet$ The case $\vert\alpha\vert=2d+2$ is covered by Lemma \ref{lemma3}.\\

$\bullet$ Consider an arbitrary $y_\alpha$ with $\vert\alpha\vert=2d$.
From the definition of the moment matrix, one may choose a pair $(i,j)$
such that the position $(i,j)$ in the matrix $M_{d+1}(y)$ lies in the submatrix
$V_{2d}$, and the corresponding entry is $y_\alpha$.
From $M_{d+1}(y)\succeq0$,
\[M_{d+1}(y)(i,i)\,M_{d+1}(y)(j,j)\,\geq\,y_\alpha^2,\]
As $M_{d+1}(y)(i,i)$ is an element $y_\beta$ of $S_{2d-2}$ with $\vert\beta\vert=2d-2$, invoking the induction hypothesis
yields $M_{d+1}(y)(i,i)\leq \tau_d^{(2d-2)/2d}$.
On the other hand, $M_{d+1}(y)(j,j)$ is a diagonal element $y_{2\beta}$ of $S_{2d+2}$ with $\vert\beta\vert=d+1$. From Lemma \ref{lemma3},
every diagonal element of $S_{2d+2}$ is dominated by  $\tau_{d+1}$, and so $M_{d+1}(y)(j,j)\leq \tau_{d+1}$. Combining the two yields
\[y_\alpha^2\,\leq\,\tau_d^{(d-1)/d}\tau_{d+1},\qquad\forall\,\alpha:\quad\va\,=\,2d.\]
Next, picking up the element $\alpha$ such that $y_\alpha=\tau_d$ one obtains
\begin{equation}
\label{tau}
\tau_d^2\,\leq\,\tau_d^{1-1/d}\tau_{d+1}\quad\Rightarrow\quad
\tau_d^{1/d}\leq\tau_{d+1}^{1/(d+1)},
\end{equation}
and so,
\[y_\alpha^2\,\leq\,\tau_d^{(d-1)/d}\tau_{d+1}\,;\qquad
\vert y_\alpha\vert^{1/\vert\alpha\vert}\,\leq\,\tau_{d+1}^{1/(2d+2)},\quad\forall\,\alpha\::\:\va=2d.\]
$\bullet$ Next, consider an arbitrary $y_\alpha$ with $\va=2d+1$.
Again, one may choose a pair $(i,j)$
such that the position $(i,j)$ in the matrix $M_{d+1}(y)$ lies in the submatrix
$V_{2d+1}$, and the corresponding entry is $y_\alpha$.
The entry $M_{d+1}(y)(i,i)$ corresponds to an element $y_{2\beta}$ of $S_{2d}$ with
$\vert\beta\vert=d$, and so, by Lemma \ref{lemma3}, 
$M_{d+1}(y)(i,i)\leq\tau_d$; similarly
the entry $M_{d+1}(y)(j,j)$ corresponds to an element $y_{2\beta}$ of $S_{2d+2}$ with
$\vert\beta\vert=d+1$, and so, by Lemma \ref{lemma3} again,
 $M_{d+1}(y)(j,j)\leq\tau_{d+1}$.
From $M_{d+1}(y)\succeq0$, we obtain
\[\tau_{d+1}\,\tau_d\,\geq\,M_{d+1}(y)(i,i)\,M_{d+1}(y)(j,j)\,\geq\,y_\alpha^2,\]
which, using (\ref{tau}), yields 
$\vert y_\alpha\vert^{1/\vert\alpha\vert}=\vert y_\alpha\vert^{1/(2d+1)}
\leq \tau_{d+1}^{1/(2d+2)}$ for all $\alpha$ with $\va=2d+1$.

$\bullet$ Finally, for an arbitrary $y_\alpha$ with $1\leq\va<2d$, use the induction hypothesis $\vert y_\alpha\vert^{1/\vert\alpha\vert}\,\leq\,\tau_d^{1/2d}$ and (\ref{tau}) to obtain
$\vert y_\alpha\vert^{1/\va}\leq\tau_{d+1}^{1/2(d+1)}$. This argument 
is also valid for the case $\vert\alpha\vert=2d$, but this latter case was treated separately to obtain (\ref{tau}).

It remains to prove that the induction hypothesis is true for $d=1$. 
This easily follows from the definition of the moment matrix $M_1(y)$.
Indeed, with $\va=1$ one has $y_\alpha^2\leq y_{2\alpha}\leq\tau_1$ (as $L_y(1)=1$), so that
$\vert y_\alpha\vert\leq \tau_1^{1/2}$ for all $\alpha$ with $\va=1$. With
$\va=2$, say with $\alpha_i+\alpha_j=2$, one has 
\[\tau_1^2\,\geq\,L_y(X_i^2)\,L_y(X_j^2)\,\geq\,L_y(X_iX_j)^2\,=\,y_\alpha^2,\]
and so $\vert y_\alpha\vert\leq \tau_1$ for all $\alpha$ with $\va=2$.
$\qed$

\subsection{Proof of Theorem \ref{th1}}
\label{proof-th1}
From (\ref{sos}), it suffices to show that $L_y(f)\geq0$ for {\it any} 
$y$ such that $M_d(y)\succeq0$, and by Remark \ref{rem1}, we may
and will assume that $L_y(1)=1$.

So let $y$ be such that $M_d(y)\succeq0$ with $L_y(1)=1$.
Let $\tau_d$ be as in Lemma \ref{lemma1} and consider the two cases $\tau_d\leq1$ and $\tau_d>1$.\\

$\bullet$ The case $\tau_d\leq 1$. By Lemma \ref{lemma1}, $\vert L_y(X^\alpha)\vert\leq 1$ for all $\alpha\in\N^n$ with
$\va \leq 2d$. Therefore,
\[L_y(f)\:\geq\quad f_0\quad-\sum_{\alpha\not\in\Gamma}
\vert f_\alpha\vert\quad+\sum_{\alpha\in\Gamma}\min[0,f_\alpha]\quad\geq0,\]
where the last inequality follows from (\ref{th1-1}).

$\bullet$ The case $\tau_d>1$. 
Recall that $L_y(1)=1$, and from Lemma \ref{lemma2}, one has
$\vert L_y(X^\alpha)\vert^{1/\va}\leq\tau_d^{1/2d}$ for all
$\alpha\in\N^n$ with $1\leq\va\leq 2d$. Therefore,
\begin{eqnarray*}
L_y(f)&\geq&f_0+(\min_{i=1,\ldots,n}f_{i2d})\,\tau_d\\
&&-\sum_{\alpha\not\in\Gamma}\vert f_\alpha\vert\,
\tau_d^{\va/2d}
+\sum_{\alpha\in\Gamma}\min[0,f_\alpha]\,\tau_d^{\va/2d}
\end{eqnarray*}
Consider the univariate polynomial $t\mapsto p(t)$, with
\[p(t)=f_0+(\min_{i=1,\ldots,n}f_{i2d})\,t^{2d}-
\sum_{\alpha\not\in\Gamma}\vert f_\alpha\vert\,t^{\vert\alpha\vert}
+\sum_{\alpha\in\Gamma}\min[0,f_\alpha]\,
t^{\vert\alpha\vert},\]
and denote $p^{(k)}\in\R[X]$, its $k$-th derivative.

By (\ref{th1-2}), $\displaystyle\min_{i=1,\ldots,n}f_{i2d}\geq0$ and so by (\ref{th1-1}), 
$p(1)\geq0$. By (\ref{th1-2}) again, $p'(1)\geq0$. In addition,
with $1\leq k\leq 2d$, (\ref{th1-2}) also implies
\begin{eqnarray*}
\min_{i=1,\ldots,n}f_{i2d}&\geq&
\sum_{\alpha\not\in\Gamma;\:\va\geq k}
\vert f_\alpha\vert
 \frac{\vert\alpha\vert}{2d}\,\frac{(\vert\alpha\vert-1)}{2d-1}\,\cdots
\frac{(\vert\alpha\vert-(k-1))}{2d-(k-1)}\\
&&-
\sum_{\alpha\in\Gamma;\:\va\geq k}
\min[0,f_\alpha]
 \frac{\vert\alpha\vert}{2d}\,\frac{(\vert\alpha\vert-1)}{2d-1}\,\cdots
\frac{(\vert\alpha\vert-(k-1))}{2d-(k-1)}
\end{eqnarray*}
because $\va-j\leq 2d-j$, for all $j=1,\ldots,k-1$, and so
\begin{eqnarray*}
\left(\prod_{j=0}^{k-1} (2d-j)\right)\min_{i=1,\ldots,n}f_{i2d}&\geq&
\sum_{\alpha\not\in\Gamma;\:\va\geq k}
\vert f_\alpha\vert\left(\prod_{j=0}^{k-1} (\va-j)\right)\\
 &-&
\sum_{\alpha\in\Gamma;\:\va\geq k}
\min[0,f_\alpha]\,
\left(\prod_{j=0}^{k-1} (\va-j)\right),
 \end{eqnarray*}
which implies $p^{(k)}(1)\geq0$.
Therefore, $p^{(k)}(1)\geq 0$ for all $k=0,1,\ldots,2d$, and so, 
as a special case of Budan-Fourier's theorem,
$p$ has no root in $(1,+\infty)$; see Basu et al \cite[Theor. 2.36]{basu}. 
Hence, $p\geq0$ on $(1,+\infty)$ and as $\tau_d>1$, 
$L_y(f)\geq p(\tau_d^{1/2d})\geq 0$. $\qed$

\subsection{Proof of Corollary \ref{coro1}}
\label{proof-coro1}
Let $y$ be such that $M_d(y)\succeq0$. Again in view of Remark \ref{rem1}, we may and will assume that $L_y(1)=1$.

Then $L_y(f)\geq\displaystyle\sum_{k=1}^d A_k$, with
\begin{eqnarray}
\label{sk}
A_k&:=&\frac{f_0}{d}+\sum_{i=1}^nf_{i2k}\,L_y(X_i^{2k})+
\sum_{\alpha\in\Gamma^2_k} \min[0,f_\alpha] \,L_y(X^\alpha) \\
\nonumber
&&-\sum_{\alpha\in\Gamma^1_k\setminus\Gamma^2_k} \vert f_\alpha\vert\, \vert L_y(X^\alpha)\vert,\qquad k=1,\ldots,d.
\end{eqnarray}
Fix $k$ arbitrary in $\{1,\ldots,d\}$ and consider the moment matrix $M_k(y)\succeq0$, which is a submatrix of $M_d(y)$. 

$\bullet$ Case $\tau_k\leq 1$. By Lemma \ref{lemma1} applied to $M_k(y)$, $\vert L_y(X^\alpha)\vert\leq 1$ for all $\alpha\in\N^n$ with
$\va \leq 2k$. Therefore, with $A_k$ as in (\ref{sk}),
\[A_k\,\geq\quad\frac{f_0}{d}\quad-\sum_{\alpha\in\Gamma^1_k\setminus\Gamma^2_k}
\vert f_\alpha\vert\quad+\sum_{\alpha\in\Gamma^2_k}\min[0,f_\alpha]\,\quad\geq0,\]
where the last inequality follows from (\ref{coro1-1}).

$\bullet$ Case $\tau_k>1$. From Lemma \ref{lemma2} applied to $M_k(y)$,
$\vert L_y(X^\alpha)\vert^{1/\va}\leq\tau_k^{1/2k}$ for all
$\alpha$ with $\va\leq 2k$. Therefore, 
$A_k\geq p_k(\tau_k^{1/2k})$, where $p_k\in\R[t]$, and
\[p_k(t)=\frac{f_0}{d}+t^{2k}\left(\,\min_{i=1,\ldots,n}f_{i2k}+\sum_{\alpha\in\Gamma^2_k}\min[0,f_\alpha]\right)- \sum_{\alpha\in\Gamma^1_k\setminus\Gamma^2_k}\vert f_\alpha\vert\,t^{\vert\alpha\vert} .\]
As in the proof of Theorem \ref{th1}, but now using (\ref{coro1-1})-(\ref{coro1-2}),
one has $p_k^{(j)}(1)\geq0$ for all 
$j=0,1,\ldots,2k$. As a particular case of Budan-Fourier's theorem,
$p_k$ has no root in $(1,+\infty)$; see Basu et al \cite[Theor. 2.36]{basu}. 
Therefore, $p_k\geq0$ on $(1,+\infty)$ which in turn implies 
$A_k\geq p_k(\tau_k^{1/2k})\geq0$ because $\tau_k>1$. Finally,
$L_y(f)\geq\sum_{k=1}^dA_k\geq0$, as $A_k\geq0$ in both cases $\tau_k\leq1$ and $\tau_k>1$. 
$\qed$

\end{document}